\documentclass[11pt]{article}

\usepackage{amsthm}
\usepackage{amsmath}
\usepackage{amssymb}
\usepackage{url}

\usepackage{graphicx}

\newtheorem{theorem}{Theorem}[subsection]
\newtheorem{lemma}[theorem]{Lemma}
\newtheorem{defn}[theorem]{Definition}
\newtheorem{prop}[theorem]{Proposition}
\newtheorem{cor}[theorem]{Corollary}

\numberwithin{equation}{section}

\newenvironment{pf}%
  {\noindent\textbf{Proof.}}%
  {\qed}

\newenvironment{pf*}%
  {\noindent\textbf{Proof.}}%

\newcommand{\Ahat}{\hat{A}}

\newcommand{\cart}{\mathcal{C}}
\newcommand{\conform}{\mathcal{G}}

\newcommand{\curvform}{\mathcal{R}}
\newcommand{\dist}{\mathcal{D}}
\newcommand{\EHcon}{\mathcal{Z}}
\newcommand{\eqn}{\mathcal{E}}
\newcommand{\eval}[2]{\left. #1 \right|_{#2}}
\newcommand{\FN}{Fr\"{o}licher--Nijenhius}
\newcommand{\fnb}[1]{[\hspace{-0.1em}[ #1 ]\hspace{-0.1em}]}
\newcommand{\biggfnb}[1]{\biggl[\hspace{-0.25em}\biggl[ #1 \biggr]\hspace{-0.25em}\biggr]}
\newcommand{\force}{\mathcal{F}}

\newcommand{\Hbar}{\bar{H}}
\newcommand{\horproj}{\mathcal{P}}
\newcommand{\identity}{\mathcal{I}}
\newcommand{\Kform}{\mathcal{K}}
\newcommand{\Lie}{\mathcal{L}}
\newcommand{\NR}{Nijenhuis--Richardson}
\newcommand{\nrb}[1]{\{\hspace{-0.3em}\{ #1 \}\hspace{-0.3em}\}}
\newcommand{\ombar}{\bar{\omega}}

\newcommand{\p}{\partial}
\newcommand{\pd}[2]{\frac{\p #1}{\p #2}}
\newcommand{\pdb}[3]{\frac{\p^2 #1}{\p #2 \p #3}}

\newcommand{\pdeform}{\mathcal{G}}
\newcommand{\Phibar}{\bar{\Phi}}

\newcommand{\Prnl}{\Pr\nolimits}
\newcommand{\psibar}{\bar{\psi}}
\newcommand{\R}{\mathbb{R}}

\newcommand{\shat}{\hat{s}}
\newcommand{\Zupw}{Z^{\scriptscriptstyle{\wedge}}\!}
\newcommand{\vertproj}{\mathcal{Q}}

\newcommand{\vf}[1]{\frac{\p}{\p #1}}
\newcommand{\vfe}[2]{\eval{\vf{#1}}{#2}}

\DeclareMathOperator{\annih}{annih}

\DeclareMathOperator{\id}{id}

\DeclareMathOperator{\Span}{span}

\newcommand{\art}[6]{#1: #2 {\it #3\/} {\bf #4} (#5) #6}
\newcommand{\book}[4]{#1: {\it #2\/} (#3, #4)}
\newcommand{\inbook}[7]{#1: #2 {\it In:\ #3\/}, ed. #4 (#5, #6) #7}
\newcommand{\arxiv}[4]{#1: #2, arXiv:#3 (#4)}

\def\im{\operatorname{Im}}
\def\v{\operatorname{v}}
\def\Tr{\operatorname{Tr}}

\begin{document}

\title{Shape maps for second order partial differential equations}

\author{O.\,Rossi, D.J.\,Saunders and G.E.\,Prince
\thanks{Research supported by grant no 14-02476S `Variations, Geometry and Physics' of the Czech Science Foundation and IRSES project GEOMECH (EU FP7, nr 246981).
Two of us (OR, DS) wish to acknowledge the hospitality of the Australian Mathematical Sciences Institute.}}

\date{}

\maketitle

\begin{abstract}
We analyse the singularity formation of congruences of solutions of systems of second order PDEs via the construction of \emph{shape maps}. The trace of such maps represents a congruence volume whose collapse we study through an appropriate evolution equation, akin to Raychaudhuri's equation. We develop the necessary geometric framework on a suitable jet space in which the shape maps appear naturally associated with certain linear connections. Explicit computations are given, along with a nontrivial example.
\end{abstract}


\section{Introduction}

This paper concerns the analysis of singularity formation of solutions of partial differential equations (PDEs). Our study is motivated by the evident relationship between that process and the collapse of dimension of congruences of geodesics on a Riemannian manifold which have been so important in cosmology \cite{DNK09,HE73,KS07,Na05,Po04,Ray55}. To get an immediate sense of what we have mind, we invite the reader now to examine Figure 1 in Section 4. The surface shown there represents the evolution of a lemniscate whose amplitude varies harmonically in time. The surface collapses both spatially and temporally, and the goal is to predict the formation of those singularities from underlying PDEs and sufficient information about the congruence.

To be more specific, we will utilise our recent work on geometric PDEs \cite{SRP15} to construct \emph{shape maps} for PDEs. The equation for the evolution of the trace of these tangent space endomorphisms is the generalization of Raychaudhuri's equation used so effectively by Hawking, Penrose and Ellis in the 1970s~\cite{HE73,HP70}. The shape map approach to Raychaudhuri's equation is due originally to Crampin and Prince \cite{CP84} in the geodesic context, and later to Jerie and Prince \cite{JePr2000} who generalized the shape map and the singularity analysis to arbitrary systems of second order ordinary differential equations (SODEs).

While our objectives are straightforward the technical tools are rather complex. So far we have restricted ourselves to PDEs which are of \emph{ connection type}, that is, second order PDEs in $m$ dependent variables and $n$ independent variables which prescribe all the second derivatives in normal form. Such equations are the natural generalizations of SODEs.

In section 2 we create the geometric framework for our PDEs. We demonstrate the interpretation of these PDEs as connections on a suitable jet space and show how they interact with the natural \emph{vertical endomorphisms} to produce various splittings of tangent spaces, identifying horizontal and vertical distributions and the corresponding projectors. In turn these allow us to define a number of curvatures and to show how they are related to the evolution of the projectors.

In section 3 we restrict attention to congruences of solutions to our PDEs and demonstrate the construction of the associated shape maps. Our first main result, Theorem~\ref{thmA}, produces an equation for the propagation of these maps in tangent directions to our solutions.  Our second main result, Theorem~\ref{TSMT}, demonstrates that the evolution of a natural measure of congruence volume depends critically on the evolution of the trace of these shape maps.

We then show how these shape maps can be combined into a vector-valued two-form, the \emph{total shape operator}, whose evolution is described in our last main result, Theorem~\ref{thmQA}. We finish off with the concrete example mentioned above.


\section{Geometric structures for second order differential equations}

In this paper we are interested in systems of overdetermined second order differential equations in several independent variables of the specific form
\begin{equation}
\label{Econnpde}
\pdb{u^\sigma}{x^i}{x^j} = F^\sigma_{ij} \biggl( x^k, u^\nu, \pd{u^\nu}{x^k} \biggr) \, , \qquad
1 \le i, j, k \le n, \; 1 \leq \sigma, \nu \leq m \, .
\end{equation}
Equations of this kind represent a direct generalization to many independent variables of systems of ordinary differential equations having the geometric meaning of second order vector fields on manifolds. Naturally, if an appropriate geometric background is used, one can discover many structures and properties of these PDEs, similar to those associated with ODEs, which remain hidden in a pure analysis approach. The aim of this paper is to develop concepts and methods, suitable for studying the behaviour of congruences of solutions of second order PDEs, which generalize the Raychaudhuri type approaches referred to in the introduction.

We now briefly describe the basics of the ODE case:\ see the book chapter~\cite{KrPr08} for details. In what follows, and in keeping with convention,  $t$ is the independent variable (rather than $x$) and $x^a$ are the dependent variables (rather than $u^a$). We start with a system of second-order ordinary differential equations
\begin{equation}
\label{SODE1}
\ddot x^a = F^a (t, x^b, \dot x^b) \, , \qquad a, b = 1, \dots, n,
\end{equation}
for some smooth functions $F^a$ on a manifold $M$ with generic local coordinates $(x^a)$. The evolution space $E:=\R \times TM$ has adapted coordinates $(t,x^a,\dot x^a)$. Using the system of equations~\eqref{SODE1} we construct on $E$ a second order differential equation field (SODE):
\begin{equation}
\label{SODE2}
\Gamma := \frac{\partial}{\partial t} + \dot x^a
\frac{\partial}{\partial x^a} + F^a \frac{\partial}{\partial \dot x^a} \, .
\end{equation}
The SODE produces on $E$ a nonlinear connection with components $\Gamma_b^a := -{\frac{1}{2}} \partial F^a / \partial \dot x^b$.

The evolution space $E$ has a number of natural structures. The contact and vertical structures are combined in the \emph{vertical endomorphism} $S$, a (1,1) tensor field on $E$, with coordinate expression $S:=V_a\otimes\omega^a$, where $V_a := \partial / \partial \dot x^a$ are the vertical basis fields and $\omega^a: = dx^a - \dot x^adt$ are local contact forms. It is shown in~\cite{CPT84} that the first order deformation $\mathcal{L}_\Gamma S$ has eigenvalues $0, 1$ and $-1$. The eigendistributions corresponding to eigenvalues $0, 1$ and $-1$ are respectively $\Span\{\Gamma\}$, the \emph{vertical distribution} $V(E)=\Span\{V_a\}$ and the \emph{horizontal distribution} $H(E)=\Span\{H_a\}$, where
\begin{equation*}
H_a := \frac{\partial}{\partial x^a} -\Gamma^b_a\frac{\partial}{\partial \dot x^b} \, .
\end{equation*}
Importantly for us there is then a canonical splitting of $TE$:
\begin{equation*}
TE = \Span\{\Gamma\} \oplus H(E) \oplus V(E) \, .
\end{equation*}


\subsection{Setting and notation}

Higher order differential equations can be represented conveniently as submanifolds of jet bundles \cite{Sau89}. We are interested in second order PDEs, so we consider a fibred manifold  $\pi:Y \to X$ where $\dim X = n >1$ and $\dim Y = n+m$, jet prolongations $\pi_1: J^1\pi \to X$ and $\pi_2: J^2\pi \to X$, and the related jet projections $\pi_{2,1}: J^2\pi \to J^1\pi$, $\pi_{1,0}: J^1\pi \to Y$ and $\pi_{2,0} = \pi_{1,0} \circ \pi_{2,1}$. For convenience we assume that all manifolds and mappings are smooth.

If $\gamma: U \to Y$ is a local section of $\pi$, so that $\pi \circ \gamma = \id_U$, then we naturally obtain sections $j^1\gamma: U \to J^1\pi$ and $j^2 \gamma: U \to J^2\pi$, called prolongations of $\gamma$.

We use local fibred coordinates. We start with coordinates $(x^i, y^\sigma)$ on an open set $V \subset Y$ adapted to the submersion $\pi$, and use associated coordinates $(x^i, y^\sigma, y^\sigma_j)$ on $\pi_{1,0}^{-1}(V) \subset J^1\pi$, and $(x^i, y^\sigma, y^\sigma_j, y^\sigma_{jk})$ on $\pi_{2,0}^{-1}(V) \subset J^2\pi$, where $1 \leq \sigma \leq m$, $1 \leq j \leq k \leq n$, and the $y^\sigma_j$ and $y^\sigma_{jk}$ arise as partial derivatives of components of local sections $\gamma$ of $\pi$:
\begin{equation*}
y^\sigma_j \circ j^1\gamma = y^\sigma_j \circ j^2\gamma = \pd{(y^\sigma \circ \gamma)}{x^j} \, ,
\qquad y^\sigma_{jk} \circ j^2\gamma = \pdb{(y^\sigma \circ \gamma)}{x^j}{x^k} \, .
\end{equation*}
Note that, for the second derivative coordinates on $J^2\pi$, the functions $y^\sigma_{ij}$ are symmetric in the indices $i$ and $j$, so that we may take those functions $y^\sigma_{jk}$ where $1 \leq j \leq k \leq n$ as independent coordinate functions. Nevertheless, when the summation convention over repeated indices is used, a summation running over \emph{all indices} will be applied.

The tangent bundle $TY$ has a vertical sub-bundle $V\pi$ containing, for every $y \in Y$, vectors  $\xi$ in $T_yY$ such that $T_y\pi (\xi) = 0$. Similarly, $TJ^1\pi$ has the vertical the sub-bundle $V\pi_1$, containing all vectors vertical with respect to the projection $\pi_1$ onto $X$, and the ``very vertical'' sub-bundle $V\pi_{1,0}$, containing vectors vertical with respect to the projection $\pi_{1,0}$ onto $Y$.

We shall need to use the natural \emph{contact structure} over the fibred manifold $\pi$, defined by the \emph{contact} (or \emph{Cartan}) \emph{distribution} annihilated by \emph{contact $1$-forms}.  On $J^1\pi$ we have the Cartan distribution of order $1$,
\begin{equation*}
\cart_{\pi_1} = \annih \{ \omega^\sigma \}
= \Span \biggl\{ \vf{x^i} + y^\sigma_i \vf{y^\sigma} \, , \;  \vf{y^\nu_k} \biggr\} \, ,
\end{equation*}
and on $J^2\pi$ the Cartan distribution of order $2$,
\begin{equation*}
\cart_{\pi_2} = \annih \{ \omega^\sigma, \; \omega^\sigma_i \}
= \Span \biggl\{ \vf{x^i} + y^\sigma_i \vf{y^\sigma} + y^\sigma_{ij} \vf{y^\sigma_{j}} \, , \; \vf{y^\nu_{kl}} \biggr\} \, ,
\end{equation*}
where
\begin{equation*}
\omega^\sigma =  dy^\sigma - y^\sigma_j dx^j \, , \quad \omega^\sigma_i =  dy^\sigma_i - y^\sigma_{ij} dx^j \, , \qquad
1 \leq \sigma \leq m \, , \;  1 \leq i \leq  n
\end{equation*}
are canonical local contact $1$-forms.

We denote by $\cart^\circ_{\pi_1} = \Span \{ \omega^\sigma \}$ and $\cart^\circ_{\pi_2} = \Span \{ \omega^\sigma, \omega^\sigma_i \}$ the annihilators of the two Cartan distributions, and more generally we write $D^\circ$ for the annihilator of a distribution $D$.

Another structure on jet bundles, fundamental for the geometry of differential equations, is the family of \emph{vertical endomorphisms}. For the first order case, each of these is an endomorphism of $TJ^1\pi$ sending, at each point $p \in J^1\pi$, vectors in $T_pJ^1\pi$ to the ``very vertical" vectors in $V_p \pi_{1,0}$, so that the endomorphism is represented by a $(1,1)$ tensor field (that is, a vector-valued $1$-form). If $\dim X = n > 1$, as in our case, each mapping arises from the choice of a non-vanishing $1$-form $\varphi$ on $X$; we then obtain an endomorphism $S^\varphi$ of $TJ^1\pi$ whose image is an $m$-dimensional sub-bundle $V^\varphi \pi_{1,0}$ of  the bundle $V\pi_{1,0} \subset TJ^1\pi$. If $\varphi = \varphi_i dx^i$ then
\begin{equation*}
S^\varphi = \varphi_i \, \omega^\sigma \otimes \vf{y^\sigma_i} \, ,
\end{equation*}
and $V^\varphi \pi_{1,0} = \Span \{ \varphi_i \, \partial / \partial y^\sigma_i \}$.

\medskip

Throughout the paper we shall frequently use vector-valued differential forms and their corresponding bracket, namely the \FN\ bracket~\cite{FrNi56}, which we denote by $\fnb{\cdot,\cdot}$. This bracket generalizes the Lie bracket of vector fields, which in this context are vector valued $0$-forms. The \FN\ bracket is given by the following formula:\ for $\omega \otimes \xi$ and $\eta \otimes \zeta$, where $\omega, \eta$ are differential forms and $\xi,\zeta$ are vector fields,
\begin{align*}
\fnb{\omega \otimes \xi, \eta \otimes \zeta}
& = \omega \wedge \eta \otimes [\xi,\zeta] + \omega \wedge \Lie_\xi \eta \otimes \zeta - \Lie_\zeta \omega \wedge \eta \otimes \xi \\
& \qquad\quad + (-1)^{\deg \omega} (d \omega \wedge i_\xi\eta \otimes \zeta + i_\zeta \omega \wedge d \eta \otimes \xi) \, .
\end{align*}
Note that if $\omega$ and $\eta$ are $1$-forms then $\fnb{\omega \otimes \xi, \eta \otimes \zeta} = \fnb{\eta \otimes \zeta, \omega \otimes \xi}$.

Now on a general manifold $M$ consider two complementary distributions of constant rank, $C$ and $D$, so that $TM = C \oplus D$. On such a manifold, every vector valued $p$-form $\Kform \in \bigwedge^p M \otimes TM$ splits into a sum $\Kform = \Kform_C + \Kform_D$, where $\Kform_C \in \bigwedge^p M \otimes C$ and $\Kform_D \in \bigwedge^p M \otimes D$ and the projectors $\Prnl_C: TM \to C$ and $\Prnl_D: TM \to D$ to the sub-bundles $C$ and $D$ extend to operators on vector valued forms such that $\Prnl_C \Kform = \Kform_C$ and $\Prnl_D \Kform = \Kform_D$.

\medskip

In what follows the \emph{connection form} of a distribution is its projector considered as a vector-valued one-form.

\begin{defn}
Given a splitting $TM = C \oplus D$ to complementary distributions on $M$, denote by $\dist$ the connection form of the distribution $D$. The vector-valued $2$-form $\curvform = \Pr_C \fnb{\dist,  \dist}$ is called \emph{curvature} of $D$ in the direction $C$. \qed
\end{defn}


\subsection{Partial differential equations of connection type}

In our standard geometric setting a system of second order PDEs~\eqref{Econnpde} is a submanifold $\eqn$ of $J^2\pi$ of dimension $n+m+nm$, defined locally by equations
\begin{equation}
\label{esub}
y^\sigma_{ij} - F^\sigma_{ij}(x^k, y^\nu, y^\nu_k) = 0 \, .
\end{equation}
The solutions of equations~\eqref{Econnpde}, if they exist, are then local sections $\gamma: U \to Y$ of the fibred manifold $\pi$ such that $j^2 \gamma(U) \subset \eqn$. It is important to stress that every section $\gamma$ of $\pi$ is, locally, a mapping $U \to U \times W \subset Y$ and so is the \emph{graph of a mapping $u: U \to W$};  in components $u^\sigma = y^\sigma \circ \gamma$, that is, $\gamma(x^i) = (x^i, u^\sigma(x^i))$.

Moreover, systems of PDEs of the form~\eqref{Econnpde} have an important geometric interpretation as \emph{jet connections}, or, equivalently, as \emph{distributions on $J^1\pi$}. In turns out that due to this property these equations are, compared to other systems of PDEs, relatively simple and manageable. Indeed, a system of overdetermined second order PDEs~\eqref{Econnpde} can be represented in different equivalent ways as follows:
\begin{enumerate}
\item As a second order Ehresmann connection, that is, a section
\begin{equation*}
\Gamma: J^1\pi \to J^2\pi \, .
\end{equation*}
In components, $y^\sigma_{ij} \circ \Gamma = F^\sigma_{ij} ( x^k, y^\nu, y^\nu_k)$, where $F^\sigma_{ij} = F^\sigma_{ji}$.
\item As a distribution $D_\Gamma$ on $J^1Y$, horizontal with respect to the projection onto $X$,
\begin{equation*}
D_\Gamma = \Span \{\Gamma_i, \, 1 \leq i \leq n\} = \annih\{ \omega^\sigma, \; \ombar_j^\sigma, \; 1 \leq \sigma \leq m, \, 1 \leq j \leq n\} \, ,
\end{equation*}
where
\begin{equation*}
\Gamma_i = \vf{x^i} + y^\sigma_i \vf{y^\sigma} + F^\sigma_{ij} \vf{y^\sigma_j}
\end{equation*}
and
\begin{equation*}
\omega^\sigma = dy^\sigma - y^\sigma_i dx^i \, \quad \ombar^\sigma_j = \Gamma^* \omega^\sigma_j = dy^\sigma_j - F^\sigma_{jk} dx^k \, ,
\end{equation*}
with $F^\sigma_{ij} = F^\sigma_{ji}$.
\item As a vector-valued form, the connection form, on $J^1Y$
\begin{equation*}
\conform = dx^i \otimes \Gamma_i
\end{equation*}
with the vector fields $\Gamma_i$ defined above.
\end{enumerate}
We can see that, for an equation of this form, the submanifold~\eqref{esub} of $J^2\pi$ is given by $\eqn = \im \Gamma = \Gamma(J^1Y)$, which is an immersed submanifold,  transversal to the fibres $\pi_{2,1}^{-1}(p)$, $p \in J^1Y$. This submanifold is endowed with the restriction of the Cartan distribution $\cart_{\pi_2}$, and we note that this restriction projects onto the distribution $D_{\Gamma}$.

A local section $\gamma$ of $\pi$ is called a \emph{solution} of the connection $\Gamma$ if $\Gamma \circ j^1\gamma = j^2\gamma$ on the domain of $\gamma$. In coordinates this is exactly a system of second order PDEs for components $u^\sigma (x^k)$, $1 \leq \sigma \leq m$, of $\gamma$, of the form~\eqref{Econnpde}. As a consequence, graphs of solutions of such equations are described by integral sections of $D_{\Gamma}$.

A connection $\Gamma$ is called \emph{integrable} if the corresponding distribution $D_{\Gamma}$ is completely integrable, satisfying Frobenius' Theorem, and this happens if and only if the curvature of the connection form $\conform$ is zero. The latter condition means that the vector-valued form $\conform$ satisfies $\fnb{\conform, \conform} = 0$:\ indeed,
\begin{equation*}
\fnb{\conform, \conform} = dx^i \wedge dx^j \otimes [\Gamma_i, \Gamma_j] \, ,
\end{equation*}
and this is zero if and only if $[\Gamma_i, \Gamma_j] = 0$ for all $i,j$. Note that, whether or not it vanishes, the Lie bracket $[\Gamma_i, \Gamma_j]$ always is a very vertical vector field, taking its values in $V\pi_{1,0}$.


\subsection{The geometry of PDEs of connection type}

We now recall some very recent results on the geometry of PDEs of connection type~\cite{SRP15}, which we shall need below to study congruences of solutions of PDEs.
These results are generalizations of well-known results for ODEs described at the beginning of this section, where a SODE $\Gamma$ defines first a splitting of $TJ^1\pi$ and then a further splitting into the three eigendistributions of $\Lie_\Gamma S$.

Given a second order connection $\Gamma: J^1\pi \to J^2\pi$, there is always a splitting of the tangent bundle $TJ^1 \pi$, defined by the complementary distributions $D_\Gamma$ and $V\pi_1$:
\begin{equation*}
TJ^1\pi = D_\Gamma \oplus V\pi_1 \, .
\end{equation*}
We can, however, construct a finer splitting, reflecting the behaviour of a particular vertical endomorphism $S^\varphi$ along the connection $\Gamma$. Such a splitting is determined by a non-vanishing \emph{closed} $1$-form $\varphi$ and a transversal direction (that is, a rank-one distribution) on the base $X$; for convenience we shall fix the scaling and instead of a direction consider a generating vector field $v$ such that $i_v \varphi = 1$. In what follows, therefore, we shall assume that the topology of $X$ admits a global structure of this kind, and that we have chosen a fixed ``normalised" pair $(\varphi, v)$ on $X$.

Given such a vector field $v$ and a non-vanishing closed $1$-form $\varphi$ on $X$, the splitting of $TJ^1\pi$ arises by considering the eigenvalue problem for the vector valued $1$-form $\Lie_{\Gamma_v} S^\varphi = \fnb{\Gamma_v, S^\varphi}$ where $\Gamma_v$ is a vector field in the distribution $D_\Gamma$ in the direction $v$, namely
\begin{equation*}
\Gamma_v = v^i \Gamma_i = v^i \biggl( \vf{x^i} + y^\sigma_i \vf{y^\sigma} + F^\sigma_{ij} \vf{y^\sigma_j} \biggr) \, .
\end{equation*}
The results can be summarised as follows.

\begin{theorem}
The vector-valued $1$-form $\Lie_{\Gamma_v} S^\varphi$ has three eigenvalues, $\lambda_1 = 0$, $\lambda_2 > 0$ and $\lambda_3 < 0$, corresponding to eigendistributions $\dist_0$ of rank $nm + n - m$, $\dist_+$ of rank $m$ and $\dist_-$ of rank $m$ as follows:
\begin{align*}
\dist_0 & = \Span \{ \Gamma_i, \, W_\sigma^p, \; 1 \leq i \leq n, \, 2 \leq p \leq n, \, 1 \leq \sigma \leq m \} \, , \\
\dist_+ & = S^\varphi(TJ^1\pi)  = \Span \biggl\{ \varphi_i \vf{y^\sigma_i}, \; 1 \leq \sigma \leq m \biggr\} \, , \\
\dist_- & = H = \Span \biggl\{H_\sigma = \vf{y^\sigma} + H_{\sigma k}^\nu \vf{y^\nu_k} \, , \ 1 \leq \sigma \leq m \biggr\} \, ,
\end{align*}
where
\begin{equation*}
W_\sigma^p = W^p_i \vf{y^\sigma_i} \, , \qquad v^i W^p_i = 0
\end{equation*}
is a local basis of the complement of $\dist_+$ in $V\pi_1$ determined by $v$, and where
\begin{equation*}
v^i (\varphi_i H_{\sigma k}^\nu + \varphi_k H_{\sigma i}^\nu) =
v^i \biggl( \varphi_j \pd{F^\nu_{ik}}{y^\sigma_j} - \delta^\nu_\sigma \pd{\varphi_k}{x^i} \biggr) \, .
\end{equation*}
This gives a threefold splitting
\begin{equation}
TJ^1\pi = \dist_0 \oplus H  \oplus \dist_+ \, .
\end{equation}

We see that $D_\Gamma \subset \dist_0$, $\dist_+ \subset V_{\pi_{1,0}}$ and $(\dist_0 \cap V\pi_{1,0}) \oplus \dist_+ = V\pi_{1,0}$, so that $D_\Gamma \subset (D_\Gamma \oplus H) \subset TJ^1\pi$ is a multiconnection (see~\cite{Sau97}) defining another threefold splitting,
\begin{equation}
\label{splmul}
TJ^1\pi = D_\Gamma \oplus H \oplus V\pi_{1,0} \, .
\end{equation}
Combining both splittings then gives a fourfold splitting
\begin{equation}
\label{spl4}
TJ^1\pi = D_\Gamma \oplus H \oplus (\dist_0 \cap V\pi_{1,0}) \oplus \dist_+
\end{equation}
into subdistributions of ranks $n$, $m$, $m$, and $nm-m$. \qed
\end{theorem}

In order to obtain explicit formulas for the eigendistributions it is convenient to take advantage of the pair $(\varphi,v)$ and use adapted coordinates $x^i$ on $X$ (and corresponding coordinates on $Y$ and $J^1\pi$) satisfying
\begin{equation*}
\varphi = dx^1 \, , \quad  i_v \varphi = v^1 = 1 \, ;
\end{equation*}
the existence of such coordinates is guaranteed by Frobenius' theorem, as $\varphi$ annihilates an integrable distribution of rank $n-1$ on $X$. We continue with the convention that indices $i, j, k, \ldots$ run from $1$ to $n$, but that an index $p$ (or $q$) runs from $2$ to $n$, both in the summation convention and when labelling families of equations. In these adapted coordinates we find that
\begin{align*}
\dist_0 & = \Span \biggl\{ \Gamma_i \, , \; W_\nu^p
= \vf{y^\nu_p} - v^p \vf{y^\nu_1} \, , \; 1 \leq i \leq n, \, 2 \leq p \leq n, \, 1 \leq \nu \leq m \biggr\} \, , \\
\dist_+ & = \Span \biggl\{ \vf{y^\sigma_1}, \; 1 \leq \sigma \leq m \biggr\} \, , \\
\dist_- & = H = \Span \biggl\{ H_\sigma = \vf{y^\sigma} + H_{\sigma k}^\nu \vf{y^\nu_k} \, , \; 1 \leq \sigma \leq m \biggr\} \, ,
\end{align*}
where
\begin{equation*}
H_{\sigma 1}^\nu = \tfrac{1}{2} \biggl( \pd{F^\nu_{11}}{y^\sigma_1} - v^p v^q \pd{F^\nu_{p q}}{y^\sigma_1} \biggr) \, ,
\quad
H_{\sigma p}^\nu = \pd{F^\nu_{1p}}{y^\sigma_1} +  v^q \pd{F^\nu_{pq}}{y^\sigma_1} = v^k \pd{F^\nu_{pk}}{y^\sigma_1} \, ;
\end{equation*}
it follows that
\begin{equation*}
2v^k H_{\sigma k}^\nu = \pd{F^\nu_{11}}{y^\sigma_1} - v^p v^q \pd{F^\nu_{p q}}{y^\sigma_1} + 2 v^p v^k \pd{F^\nu_{pk}}{y^\sigma_1}
= v^i v^k \pd{F^\nu_{ik}}{y^\sigma_1} \, .
\end{equation*}

\medskip
In all three splittings, a central role is played by the distribution $H$; we call it a \emph{semi-horizontal} bundle associated to $\Gamma$. If now take the annihilator of $D_\Gamma \oplus H$ we obtain an important codistribution, the \emph{force bundle},
\begin{equation*}
\force = D_\Gamma^\circ \cap H^\circ = \Span \{ \psi^\nu_k, \;  1 \leq k \leq n, \, 1 \leq \nu \leq m\}
\end{equation*}
where, in adapted coordinates,
\begin{equation*}
\psi^\nu_k := dy^\nu_k - F^\nu_{ki} dx^i - H^\nu_{\sigma k} \omega^\sigma \, .
\end{equation*}

\medskip

The splittings of the tangent bundle induce dual splittings of the cotangent bundle. The threefold splitting $D_\Gamma \oplus H \oplus V\pi_{1,0}$ from~\eqref{splmul} gives a dual splitting
\begin{equation}
\label{dsplmul}
T^*J^1 \pi = \pi_1^* (T^*X) \oplus \cart^\circ_{\pi_1} \oplus \force \, ,
\end{equation}
and induces a decomposition of $1$-forms into horizontal, contact and force forms. We therefore obtain three vector valued $1$-forms giving a decomposition of the identity
\begin{equation*}
\identity = \conform + \horproj + \vertproj \, , \qquad
\conform : TJ^1\pi \to D_\Gamma \, , \quad \horproj : TJ^1\pi \to H \, , \quad \vertproj : TJ^1\pi \to V\pi_{1,0} \, .
\end{equation*}
We see immediately that the canonical local bases of vector fields and of $1$-forms, adapted to the dual splittings~\eqref{splmul} and~\eqref{dsplmul}, are
\begin{equation*}
\biggl(\Gamma_i, H_\sigma, \vf{y^\nu_k} \Bigr) \, , \quad (dx^i, \omega^\sigma, \psi^\nu_k) \, .
\end{equation*}
In these bases,
\begin{equation*}
\conform = dx^i \otimes \Gamma_i \, , \quad \horproj = \omega^\sigma \otimes H_\sigma \, ,
\quad \vertproj = \psi^\sigma_i \otimes \vf{y^\sigma_i} \, .
\end{equation*}
As we shall see, the vector valued $1$-form $\vertproj$, taking its values in the very vertical bundle $V\pi_{1,0}$, will be particularly important.

Finally, we note that the threefold splitting~\eqref{dsplmul} may be refined to give a fourfold splitting of $T^* J^1\pi$ dual to the fourfold splitting~\eqref{spl4} of $TJ^1\pi$. This may be obtained simply by decomposing the force bundle $\force$ and takes the form
\begin{equation*}
T^*J^1\pi = \pi_1^* (T^*X) \oplus \cart^\circ_{\pi_1} \oplus (\dist^\circ_+ \cap \force) \oplus \force_+ \, ,
\end{equation*}
where we have denoted by $\force_+ = \dist_0^\circ \cap \dist_-^\circ$ the $m$-dimensional sub-bundle of the force bundle complementary to $\dist^\circ_+$. The corresponding adapted basis and dual basis are
\begin{equation*}
\biggl(\Gamma_i, \, H_\sigma, \,  W_\sigma^p, \, \vf{y^\nu_1} \biggr) \, , \quad
(dx^i \, , \, \omega^\sigma \, , \, \psi^\nu_p \, , \, v^k\psi^\nu_k) \, .
\end{equation*}
Thus the projectors $\conform$ and $\horproj$ to $D_\Gamma$ and $H$ are the same, and the projector $\vertproj$ to $V\pi_{1,0}$ splits into a sum of two projectors, of which the projector $\vertproj_{+}$ to $\dist_+$ will be used below for construction of a shape operator. So now we have a decomposition of the identity
\begin{equation*}
\identity = \conform + \horproj + \tilde{\vertproj}  + \vertproj_+ \ ,
\end{equation*}
and in adapted coordinates
\begin{equation*}
\tilde{\vertproj} = \psi^\nu_p \otimes W_\nu^p \, , \quad \vertproj_+ = v^k \psi^\sigma_k \otimes \vf{y^\sigma_1} \, .
\end{equation*}


\subsection{Curvatures induced by second order PDEs}

Whenever we have a system of PDEs represented by a connection $\Gamma$, we may obtain a number of different curvature operators as explained in~\cite{SRP15}.

First we have the canonical curvature, arising from the splitting $TJ^1\pi = D_\Gamma \oplus V\pi_1$. This is a vector valued $2$-form, expressed in terms of the Lie bracket of the vector fields belonging to the distribution $D_\Gamma$ as
\begin{equation*}
{\curvform}^{\Gamma} = \Prnl_{V\pi_{1}} \fnb{\conform, \conform} = \fnb{\conform, \conform} = dx^i \wedge dx^k \otimes [\Gamma_i, \Gamma_k] \, .
\end{equation*}

There are, however, other curvature operators, obtained from the individual components of more refined splittings of $TJ^1\pi$ associated with $\Gamma$. In particular, we shall be interested in curvatures obtained from the splittings~\eqref{splmul} and~\eqref{spl4} considered in the previous section.

The threefold splitting~\eqref{splmul} gives the following curvatures.

\begin{theorem}
The splitting $TJ^1\pi = (D_\Gamma \oplus H) \oplus V\pi_{1,0}$ gives rise to the following curvature operators on $J^1\pi$.
\begin{enumerate}
\item The inner curvature of the connection $\Gamma$. This is equal to the canonical curvature, and is given by
\begin{equation*}
R^{\Gamma} = \Prnl_{V{\pi_{1,0}}} \fnb{\conform, \conform} = \fnb{\conform, \conform}
= dx^i \wedge dx^k \otimes [\Gamma_i, \Gamma_k] = {\curvform}^\Gamma \, .
\end{equation*}
\item The curvature of $\Gamma$ with respect to $H$. This is known as the Jacobi endomorphism, by analogy with the ODE case, and is given by
\begin{equation*}
\Phi = \Prnl_{V\pi_{1,0}} \fnb{\conform, \horproj} = \fnb{\conform, \horproj} - dx^i \wedge \psi^\sigma_i \otimes H_\sigma
= \Phi^\nu_{i \sigma k} dx^i \wedge \omega^\sigma \otimes \vf{y^\nu_k}
\end{equation*}
where
\begin{equation*}
\Phi^\nu_{i \sigma k} = \Gamma_i(H^\nu_{\sigma k}) - H_\sigma(F^\nu_{ik}) + H^\rho_{\sigma i} H_{\rho k}^\nu \, . \tag*{\qed}
\end{equation*}
\end{enumerate}
\end{theorem}

Using $\identity = \conform + \horproj + \vertproj$ we see that $\fnb{\conform, \vertproj} = - \fnb{\conform, \conform} - \fnb{\conform, \horproj}$, and so projecting onto the very vertical bundle $V\pi_{1,0}$ we obtain the \emph{vertical curvature} of the equations.

\begin{theorem}
\begin{equation*}
\Prnl_{V{\pi_{1,0}}} \fnb{\conform, \vertproj} = \vertproj \circ \fnb{\conform, \vertproj} = - R^\Gamma - \Phi \, . \tag*{\qed}
\end{equation*}
\end{theorem}

The fourfold splitting~\eqref{spl4} gives refined versions of these operators.

\begin{theorem}
The splitting $TJ^1\pi = (D_\Gamma \oplus H \oplus (\dist_0 \cap V\pi_{1,0})) \oplus \dist_+$ gives rise to the following curvature operators on $J^1\pi$.
\begin{enumerate}
\item The inner curvature of the connection $\Gamma$ in the direction of $\dist_+$. This is a partial canonical curvature, given by
\begin{equation}
R^\Gamma_+ = \Prnl_{\dist_+} \fnb{\conform, \conform} = \Prnl_{\dist_+} R^\Gamma
= dx^i \wedge dx^j \otimes \bigl( v^k [\Gamma_i, \Gamma_j]^\sigma_k \bigr) \vf{y^\sigma_1} \, ,
\end{equation}
where, for each pair of indices $(i,j)$ the Lie bracket $[\Gamma_i, \Gamma_j]$ is a very vertical vector field whose components are $[\Gamma_i, \Gamma_j]^\sigma_k$.
\item The curvature of $\Gamma$ with respect to $H$ in the direction of $\dist_{+}$. This is a partial Jacobi endomorphism, given by
\begin{equation*}
\Phi_+ = \Prnl_{\dist_+} \fnb{\conform, \horproj}
= \Prnl_{\dist_+} \Phi = \Phi_{i \nu}^\sigma dx^i \wedge \omega^\nu \otimes \vf{y^\sigma_1}
\end{equation*}
where
\begin{equation*}
\Phi_{i \nu}^\sigma =  v^k \Phi^\nu_{i \sigma k}
= v^k \bigl( \Gamma_i(H^\sigma_{\nu k}) - H_\nu(F^\sigma_{ik}) + H^\rho_{\nu i} H_{\rho k}^\sigma \bigr) \, .
\end{equation*}
\item The curvature of $\Gamma$ with respect to $\dist_0 \cap V\pi_{1,0} = \Span \{W_\sigma^p \}$, given by
\begin{align*}
r_+ & = \Prnl_{\dist_+} \fnb{\conform, \tilde{\vertproj}} \\
& = dx^i \wedge \psi^\nu_p \otimes v^k \biggl(v^p \pd{F^\sigma_{ik}}{y^\nu_1} - \pd{F^\sigma_{ik}}{y^\nu_p}
- (v^p \delta^1_i - \delta^p_i) H_{\nu k}^\sigma \biggr) \vf{y^\sigma_1} - dx^i \wedge \psi^\nu_p \otimes \pd{v^p}{x^i} \vf{y^\nu_1} \, . \tag*{\qed}
\end{align*}
\end{enumerate}
\end{theorem}

The vertical curvature with respect to the complementary vertical bundle of the splitting is now a curvature with respect to $\dist_{+}$, and takes the following form.

\begin{theorem}
\label{thmQ}
\begin{equation*}
\Prnl_{\dist_+} \fnb{\conform, \vertproj_+} = \vertproj_+ \circ \fnb{\conform, \vertproj_+} = - R^\Gamma_+ - \Phi_+ - r_+ \, . \tag*{\qed}
\end{equation*}
\end{theorem}

We note, finally, a result concerning the evolution of the projector $\vertproj_+$.

\begin{prop}
\label{Pvpev}
The evolution of $\vertproj_+$ in the direction $\Gamma_v$ is given by
\begin{equation*}
\Lie_{\Gamma_v} \vertproj_+ = \fnb{\Gamma_v, \vertproj_+} =
i_{\Gamma_v} \fnb{\pdeform, \vertproj_+} + \tfrac{1}{2} i_{\Gamma_v} R^\Gamma_+ \, .  \tag*{\qed}
\end{equation*}
\end{prop}


\section{Shape maps for second order PDEs}

In this section we shall consider how to construct `shape maps' for systems of PDEs represented by second order connections. We shall, in fact, construct such maps on two different levels. The first type of shape map will be a vector-valued $1$-form, like the shape map associated with a second-order ODE, and will depend upon the choice of a nondegenerate pair $(\varphi,v)$ as used in the tangent bundle decompositions described above. The second will be a vector-valued $2$-form from which the vector-valued $1$-form can be obtained by contraction with $v$. Although the latter might seem to be the more general object, its construction still depends upon the choice of $v$ as well as $\varphi$, so it contains essentially the same information.

Both types of shape map are associated with particular choices of congruences of solutions, and may be used to analyse the way in which such congruences might collapse.


\subsection{Congruences of solutions}

We start by defining what we mean by a `congruence of solutions' of a family of PDEs.

Consider the system of overdetermined second order PDEs given in equations~\eqref{Econnpde}, represented by a second order connection $\Gamma: J^1\pi \to J^2\pi$. Let $Z : V \to J^1\pi$, where $V \subset $ is an open subset, be a local section; the map $Z$ is called a local connection (Ehresmann connection, or jet connection) on $\pi$, and the image submanifold $\im Z = Z(V) \subset J^1 \pi$ is an immersed submanifold transversal to the fibres $\pi_{1,0}^{-1}(y)$ for $y \in V$. We shall denote by $D_Z$ the distribution on $V$ associated with $Z$, and by $\EHcon$ the corresponding connection form. If $y^\sigma_{i} \circ Z = Z^\sigma_i(x^k, y^\nu)$ are the components of $Z$ then
\begin{equation*}
D_Z = \Span \{Z_i,  \, 1 \leq i \leq n\} = \annih\{\ombar^\sigma, \, 1 \leq \sigma \leq m\} \, ,
\quad \EHcon = dx^i \otimes Z_i
\end{equation*}
where
\begin{equation*}
Z_i = \vf{x^i} + Z^\sigma_i \vf{y^\sigma}
\end{equation*}
and where $\ombar^\sigma = Z^* \omega^\sigma = dy^\sigma - Z^\sigma_j dx^j$; we shall often use a bar to indicate the pullback by $Z$ of a function or form defined locally on $J^1\pi$. A local connection $Z$ on $\pi$ therefore represents a system of first order PDEs `in normal form',
\begin{equation*}
\pd{u^\sigma}{x^i} = Z^\sigma_{i} (x^k, u^\nu) \, ,
\end{equation*}
for the components $u^\sigma(x^k)$, $1 \leq \sigma \leq m$, of local sections of the fibred manifold $\pi$ taking their values in $V$.

A local connection on $\pi$ induces a splitting of every vector field on $V$ into a horizontal component belonging to $D_Z$ and vertical component belonging to $V\pi$, and lifts vector fields on $U$ to vector fields on the manifold $\im Z \subset J^1Y$. We note that for any such vector field $\xi$
\begin{equation}
\label{tanZ}
TZ(\xi) = \xi + \xi(Z^\nu_k) \vf{y^\nu_k} \, .
\end{equation}

Given a second order connection $\Gamma$, a local connection $Z$ defined on $V \subset Y$ is said to be \emph{embedded in $\Gamma$} if at every point $y \in V$
\begin{equation} \label{embZ}
T_y Z(D_Z(y)) \subset D_\Gamma (Z(y)) \, .
\end{equation}
This condition says that the tangent map to $Z$ transfers the distribution $D_Z$ to $D_\Gamma$ along the submanifold $\im Z \subset J^1\pi$. In terms of solutions of the corresponding PDEs, an embedded connection lifts solutions of $Z$ to solutions of $\Gamma$:\ if $\gamma$ is a solution of $Z$, that is $Z \circ \gamma = j^1 \gamma$, then $\Gamma \circ Z \circ \gamma = j^2 \gamma$.

\begin{defn} By a \emph{congruence} of (local) solutions of a system of PDEs defined by a second order connection $\Gamma$ we shall mean a (local) connection on $\pi$ embedded in $\Gamma$. \qed
\end{defn}

We stress that, in this setting, solutions of PDEs are sections of the fibred manifold $\pi$ rather than functions. This means that we are dealing with the graphs of maps $(x^i) \mapsto (u^\sigma(x^i))$ rather than with the maps themselves.

From equation~\eqref{tanZ} we can see that since the lifted distribution $D_Z = \Span \{Z_i \}$ is the restriction of the distribution $D_\Gamma = \Span \{\Gamma_i\}$ to the points of the submanifold $\im Z$, it follows that
\begin{equation*}
F^\sigma_{ij} \circ Z = Z_i (Z^\sigma_j) \, .
\end{equation*}
But the functions $F^\sigma_{ij}$ are symmetric in $i,j$, so that $Z_i (Z^\sigma_j) - Z_j (Z^\sigma_i) = [Z_i, Z_j] = 0$, in other words that $\fnb{\mathcal{Z}, \mathcal{Z}} = 0$, and so we obtain the following result.

\begin{prop}
Every embedded connection is completely integrable. \qed
\end{prop}

This means that a congruence of solutions is an $n$-dimensional foliation of the manifold $\im Z \subset J^1\pi$.


\subsection{Directional shape operators}

We shall now construct a `directional' shape operator for a congruence of solutions of a system of overdetermined second order PDEs, where the PDEs are represented by a second order connection $\Gamma$ and the congruence is represented by a first order connection $Z$ embedded in $\Gamma$. We shall assume that we have chosen a non-vanishing closed form $\varphi$ and a vector field $v$ on $X$ such that $i_v \varphi = 1$; these will be fixed until further notice, so that we have the four-fold splitting~\eqref{spl4}
\begin{equation*}
TJ^1\pi = D_\Gamma \oplus H \oplus (\dist_0 \cap V\pi_{1,0}) \oplus \dist_+ \, .
\end{equation*}
Recall that $\dist_+ = S^\varphi(V\pi_1)$; we shall use the inverse $s_{j^1_x\gamma}$ of the vertical lift map
\begin{equation*}
l_{j^1_x\gamma, \varphi} : V_{\gamma(x)}\pi \to V_{j^1_x\gamma}\pi_{1,0}
\end{equation*}
used to construct $S^\varphi$, so that $s : \dist_+ \to V\pi$ is a fibrewise isomorphism depending on $\varphi$. This will enable us to introduce a map $Z^\varphi$ from vector-valued forms on $\im(Z) \subset J^1\pi$ with values in $\dist_+$ to vertical-valued forms on $Y$. If $\mathcal{K} \in \Lambda^q(J^1\pi) \otimes \dist_+$ and $\xi_1, \xi_2, \ldots, \xi_q$ are vector fields on $Y$ then, for any $y \in Y$, put
\begin{equation*}
(Z^\varphi \mathcal{K})(y) \, (\xi_1, \xi_2, \ldots, \xi_q)
= s_{Z(y)} \Bigl( \mathcal{K} \bigl( Z(y) \bigr) \, \bigl( TZ(\xi_1), TZ(\xi_2), \ldots, TZ (\xi_q) \bigr) \Bigr) \, ,
\end{equation*}
so that $Z^\varphi \mathcal{K} \in \Lambda^q Y \otimes V\pi$. We may express this construction in coordinates adapted to $\varphi$ with $\varphi = dx^1$, so that
\begin{equation*}
\varphi_1 = 1 \, , \quad  \varphi _2 = \dots = \varphi_n = 0 \, , \qquad \text{and} \quad v^1 = 1 \, ;
\end{equation*}
then
\begin{equation*}
s_{Z(y)} \biggl( \vfe{y^\sigma_1}{Z(y)} \biggr) = \vfe{y^\sigma}{y} \, ,
\end{equation*}
and if $\mathcal{K} = \eta^\sigma \otimes \p / \p y^\sigma_1$ then
\begin{equation*}
Z^\varphi \mathcal{K} = Z^* \eta^\sigma \otimes \vf{y^\sigma} \, .
\end{equation*}

We can now define the shape operator.
\begin{defn}
The \emph{shape operator} for the congruence defined by $Z$ in the fixed direction $(\varphi, v)$ is the vector valued $1$-form $A_Z$ on $Y$ given by
\begin{equation*}
A_Z = Z^\varphi \vertproj_+ \, . \tag*{\qed}
\end{equation*}
\end{defn}
Note that, although $Z^\varphi$ does not depend on $v$, the shape operator $A_Z$ depends on both $\varphi$ and $v$, because $\vertproj$ does.

\begin{prop}
For every vector field $\xi$ on $Y$ the shape operator $A_Z$ satisfies the identity
\begin{equation*}
S^{\varphi} \bigl( TZ (A_Z(\xi)) \bigr) = \vertproj_+ (TZ(\xi))
\end{equation*}
along $\im Z$.
\end{prop}

\begin{pf*}
We shall prove that at each $p \in \im Z$, where $p = Z(y)$ for $y = \pi_{1,0}(p) \in Y$, both sides are equal to $ s_y^{-1} \bigl( A_Z (y) \, (\xi) \bigr)$.

Put
\begin{equation*}
\Xi = A_Z(\xi) \, , \qquad \Xi = \Xi^\sigma \vf{y^\sigma}
\end{equation*}
so that
\begin{equation*}
\Xi(y) = A_Z(y) \, (\xi) = (Z^\varphi \vertproj_+) (y) \, (\xi) = s_y \Bigl( \vertproj_+ \bigl( Z(y) \bigr) \, \bigl( TZ(\xi) \bigr) \Bigr) \, ;
\end{equation*}
it follows that the right-hand side is
\begin{equation*}
\vertproj_+ (Z(y)) \, (TZ(\xi)) = s_y^{-1} \bigl( A_Z (y) \, (\xi) \bigr) \, .
\end{equation*}

On the other hand,
\begin{equation*}
TZ(\Xi) = \Xi + \Xi(Z^\nu_k) \vf{y^\nu_k} \, ,
\end{equation*}
so that the left-hand side is
\begin{equation*}
S^{\varphi} \bigl( TZ (A_Z(\xi)) \bigr) = S^\varphi (TZ(\Xi)) = \omega^\sigma (\Xi) \,  \vf{y^\sigma_1} = \Xi^\sigma \vf{y^\sigma_1}
= s^{-1} (\Xi) =  s^{-1} (A_Z (\xi)) \, . \tag*{\qed}
\end{equation*}
\end{pf*}

We can now see how the shape operator acts on vector fields on $Y$.
\begin{prop}
If $\xi$ is a vector field on $Y$ then $A_Z(\xi) = A_Z(\xi_{\v})$, where $\xi_{\v} = \xi - \mathcal{Z}(\xi)$ is the vertical projection. In particular, if $\xi \in D_Z$ then $A_Z(\xi) = 0$.
\end{prop}
\begin{pf}
As $Z$ is a connection on $Y$ it defines a splitting $TY = D_Z \oplus V\pi$, so that $\xi$ can be written as the sum of a `horizontal' and vertical vector,
\begin{equation*}
\xi = \xi_{Z} + \xi_{\v} \, ,
\end{equation*}
where $\xi_Z = \EHcon(\xi)$ is the projection to $D_Z$ and $\xi_{\v}$ is the complementary vector field taking values in $V\pi(y)$. It follows that
\begin{equation*}
TZ(\xi) = TZ(\xi_Z + \xi_{\v}) = TZ(\xi_Z) + TZ(\xi_{\v}) \, .
\end{equation*}
But $TZ(\xi_Z) \in \eval{D_\Gamma}{\im Z}$, so that $\vertproj_+ (TZ(\xi_Z)) = 0$ and therefore $A_Z(\xi_Z) = 0$; thus $A_Z(\xi) = A_Z(\xi_{\v})$.
\end{pf}

\medskip
We shall need to use the coordinate expression of the shape operator, which is derived easily from its definition. We find that
\begin{equation*}
A_Z = v^i \psibar^\sigma_i \otimes \vf{y^\sigma}
= v^i \biggl( \pd{Z^\sigma_i}{y^\nu} - \Hbar^\sigma_{\nu i} \biggr) \ombar^\nu \otimes \vf{y^\sigma} \, ,
\end{equation*}
where $\psibar^\sigma_i = Z^* \psi^\sigma_i$, $\ombar^\nu = Z^* \omega^\nu$ and $\Hbar^\sigma_{\nu i} = H^\sigma_{\nu i} \circ Z$. Thus, denoting the components of $A_Z$ by $A^\sigma_{\nu}$, we have
\begin{equation*}
A^\sigma_{\nu} = v^i \biggl( \pd{Z^\sigma_i}{y^\nu} - \Hbar^\sigma_{\nu i} \biggr) \, .
\end{equation*}

We now want to discover how the shape operator $A_Z$ evolves along a congruence in the direction determined by $v$. To do this, we first investigate the evolution of the projection $\vertproj_+$.
\begin{theorem}
\label{thmQ+}
The evolution of the projection operator $\vertproj_+$ satisfies
\begin{equation*}
\vertproj_+ \circ \Lie_{\Gamma_v} \vertproj_+  = \Prnl_{V_+} \fnb{\Gamma_v, \vertproj_+} =  - i_{\Gamma_v} ( \tfrac{1}{2} R^\Gamma_+ + \Phi_+  + r_+) \,.
\end{equation*}
\end{theorem}
\begin{pf*}
From Proposition~\ref{Pvpev} and Theorem~\ref{thmQ} we have
\begin{align*}
\vertproj_+ \circ \Lie_{\Gamma_v} \vertproj_+
& = \Prnl_{V_+} \fnb{\Gamma_v, \vertproj_+} \\
& = \Prnl_{V_+} \bigl( i_{\Gamma_v} \fnb{\pdeform, \vertproj_+} + \tfrac{1}{2} i_{\Gamma_v} R^\Gamma_+ \bigr) \\
& = - i_{\Gamma_v} ( \tfrac{1}{2} R^\Gamma_+ + \Phi_+  + r_+) \, . \tag*{\qed}
\end{align*}
\end{pf*}

This is a generalization of the corresponding evolution formula for ODEs, where now on the right-hand side we have three curvature components rather than one. In the ODE case the curvature $r_+$ does not exist because $\dist_+$ is the whole of $V\pi_{1,0}$, and the curvature $R^\Gamma$ vanishes because the equation is a single vector field generating a one-dimensional distribution which is always completely integrable, so we are left with the Jacobi endomorphism.

We now let $Z_v$ denote the vector field $i_v Z$, so that
\begin{equation*}
Z_v = v^i \biggl( \vf{x^i} + Z^\sigma_i \vf{y^\sigma} \biggr) = v^i Z_i \, ;
\end{equation*}
we can see how the shape operator $A_Z$ propagates along $Z_v$.
\begin{theorem} \emph{(First Main Theorem)}
\label{thmA}

The directional shape operator $A_Z$ satisfies the equation
\begin{equation*}
\Lie_{Z_v} A_Z  = - (A_Z)^2  - i_{Z_v} (Z^\varphi \Phi_+   +  Z^\varphi r_+ ) \, .
\end{equation*}
\end{theorem}
\begin{pf}
We carry out a direct calculation. We see first that
\begin{align*}
\Lie_{Z_v} A_Z & = \fnb{Z_v, A_Z} \\
& = v^i \pd{v^k}{x^i} \psibar^\sigma_k \otimes \vf{y^\sigma} + v^i v^k \biggl\{ Z_i \biggl( \pd{Z^\sigma_k}{y^\nu} \biggr) - Z_i(\Hbar^\sigma_{\nu k}) \\
& \qquad + \biggl(\pd{Z^\sigma_k}{y^\rho} - \Hbar^\sigma_{\rho k} \biggr) \pd{Z^\rho_i}{y^\nu} - \biggl(\pd{Z^\rho_k}{y^\nu}
- \Hbar^\rho_{\nu k} \biggr) \pd{Z^\sigma_i}{y^\rho} \biggr\} \ombar^\nu \otimes \vf{y^\sigma} \, ,
\end{align*}
and that
\begin{align*}
A_Z^2 & = A^\sigma_\rho A^\rho_\nu \, \ombar^\nu \otimes \vf{y^\sigma} \\
& = v^i v^k \biggl( \pd{Z^\sigma_i}{y^\rho} \pd{Z^\rho_k}{y^\nu} - \Hbar^\sigma_{\rho i} \pd{Z^\rho_k}{y^\nu}
- \Hbar^\rho_{\nu k} \pd{Z^\sigma_i}{y^\rho} + \Hbar^\sigma_{\rho i} \Hbar^\rho_{\nu k} \biggr) \, \bar \omega^\nu \otimes \vf{y^\sigma} \, .
\end{align*}
Next we shall use the fact that $Z^* df = d(f \circ Z)$ for any function $f$ on $J^1Y$, giving
\begin{align*}
Z_i (f \circ Z) & =  \Gamma_i(f) \circ Z \, , \\
\pd{(f \circ Z)}{y^\nu} & =  H_\nu (f) \circ Z
+ \biggl( \pd{f}{y^\sigma_j} \circ Z \biggr) \biggl(\pd{Z^\sigma_j}{y^\nu} - \Hbar^\sigma_{\nu j} \biggr) \, .
\end{align*}
We can now calculate $i_{Z_v} Z^\varphi \Phi_+$. We obtain
\begin{align*}
\Phibar_{i \nu}^\sigma = \Phi_{i \nu}^\sigma \circ Z & =
v^k \Gamma_i(H^\sigma_{\nu k}) \circ Z - v^k H_\nu(F^\sigma_{ik}) \circ Z + v^k \Hbar^\rho_{\nu i} \Hbar_{\rho k}^\sigma \\
& = v^k \biggl( Z_i(\Hbar^\sigma_{\nu k}) - \pd{\bar F^\sigma_{ik}}{y^\nu}
+ \biggl( \pd{F^\sigma_{ik}}{y^\rho_j} \circ Z \biggr) \biggl( \pd{Z^\rho_j}{y^\nu} - \Hbar^\rho_{\nu j} \biggr)
+ \Hbar^\rho_{\nu i} \Hbar_{\rho k}^\sigma \biggr) \, ,
\end{align*}
so that
\begin{align*}
i_{Z_v} Z^\varphi \Phi_+ & = v^i \Phibar_{i \nu}^\sigma \, \ombar^\nu \otimes \vf{y^\sigma} \\
& = v^i v^k \biggl( Z_i(\Hbar^\sigma_{\nu k}) - \pd{\bar F^\sigma_{ik}}{y^\nu}
+ \biggl( \pd{F^\sigma_{ik}}{y^\rho_j} \circ Z \biggr) \biggl( \pd{Z^\rho_j}{y^\nu} - \Hbar^\rho_{\nu j} \biggr)
+ \Hbar^\rho_{\nu i} \Hbar_{\rho k}^\sigma \biggr) \ombar^\nu \otimes \vf{y^\sigma} \\
& = v^i v^k \biggl( Z_i(\bar H^\sigma_{\nu k} ) - \pd{Z_i(Z^\sigma_k)}{y^\nu} + \Hbar^\rho_{\nu i} \Hbar_{\rho k}^\sigma \biggr) \,
\ombar^\nu \otimes \vf{y^\sigma} \\
& \qquad\qquad + v^i v^k \biggl( \pd{F^\sigma_{ik}}{y^\rho_j} \circ Z \biggr) \, \psibar^\rho_j \otimes \vf{y^\sigma} \, .
\end{align*}
Finally we see that
\begin{align*}
i_{Z_v} Z^\varphi r_+
& = v^i v^k \psibar^\nu_p \otimes \biggl\{ \biggl(v^p \pd{F^\sigma_{ik}}{y^\nu_1} - \pd{F^\sigma_{ik}}{y^\nu_p}
- (v^p \delta^1_i - \delta^p_i) H_{\nu k}^\sigma \biggr) \circ Z \biggr\} \vf{y^\sigma} \\
& \qquad\qquad - v^i \psibar^\nu_p \otimes \biggl( \pd{v^p}{x^i} \vf{y^\nu} \biggr) \\[2ex]
& = \biggl( v^i v^k v^p \pd{F^\sigma_{ik}}{y^\nu_1} - v^i v^k \pd{F^\sigma_{ik}}{y^\nu_p}
- (v^1 v^k v^p H^\sigma_{\nu k} - v^p v^k H_{\nu k}^\sigma \biggr) \biggr) \circ Z \, \psibar^\nu_p \otimes \vf{y^\sigma} \\
& \qquad\qquad - v^i \pd{v^p}{x^i} \psibar^\nu_p \otimes \vf{y^\nu} \\[2ex]
& = v^i v^k \biggl\{ \biggl( v^p \pd{F^\sigma_{ik}}{y^\nu_1}
- \pd{F^\sigma_{ik}}{y^\nu_p} \biggr) \circ Z \biggr\} \, \psibar^\nu_p \otimes \vf{y^\sigma} - v^i \pd{v^p}{x^i} \psibar^\nu_p \otimes \vf{y^\nu}
\end{align*}
so that, adding,
\begin{align*}
\lefteqn{\Lie_{Z_v} A_Z + (A_Z)^2 + i_{Z_v} Z^\varphi \Phi_+  + i_{Z_v} Z^\varphi r_+} \\[1ex]
& \qquad = -2 v^i v^k \Hbar^\sigma_{\rho k} \, \psibar^\rho_i \otimes \vf{y^\sigma}
+ v^i v^k \biggl( \pd{F^\sigma_{ik}}{y^\rho_1} \circ Z \biggr) \, \psibar^\rho_1 \otimes \vf{y^\sigma} \\
& \qquad\qquad\qquad + v^i v^k v^p \biggl( \pd{F^\sigma_{ik}}{y^\nu_1} \circ Z \biggr) \, \psibar^\nu_p \otimes \vf{y^\sigma} \\[2ex]
& \qquad = \biggl\{ \biggr( v^i v^k \pd{F^\sigma_{ik}}{y^\rho_1}
- 2 v^k H^\sigma_{\rho k} \biggr) \circ Z \biggr\} \psibar^\rho_1  \otimes \vf{y^\sigma} \\
& \qquad\qquad\qquad +  v^p \biggl\{ \biggl( v^i v^k \pd{F^\sigma_{ik}}{y^\nu_1}
- 2 v^k  H^\sigma_{\nu k} \biggr) \circ Z \biggr\} \psibar^\nu_p \otimes \vf{y^\sigma} = 0
\end{align*}
as required.
\end{pf}


\subsection{The collapse of a congruence}

In order to test whether a congruence of solutions `collapses' we shall consider characterising $m$-forms for the distribution $D_Z$, and examine how they evolve under the flow of the vector field $Z_v$.
\begin{defn}
A characterising $m$-form $\tilde{\Omega}$ for the distribution $D_Z$ is called a \emph{$\Gamma$-compatible volume} if
\begin{equation*}
\Lie_{Z_v} \tilde{\Omega} = i_{A_Z} \tilde{\Omega} \, . \tag*{\qed}
\end{equation*}
\end{defn}
To establish some properties of $\Gamma$-compatible volumes, we start with the characterising form for $D_Z$ corresponding to a choice of coordinates.
\begin{lemma}
\label{Lshapevol}
If
\begin{equation*}
\Omega = \ombar^1 \wedge \ombar^2 \wedge \cdots \wedge \ombar^m
\end{equation*}
is the natural characterising form for the distribution $D_Z$ in the coordinates $(x^i, y^\sigma)$ on $Y$ then
\begin{equation*}
i_{A_Z} \Omega = (\Tr A_Z) \, \Omega
\end{equation*}
where $\Tr A_Z$ is the trace of the matrix of the components of $A_Z$.
\end{lemma}
\begin{pf*}
A straightforward calculation yields
\begin{equation*}
i_{A_Z} \Omega =  A^\sigma_{\nu} \ombar^\nu \wedge  i_{\p / \p y^\sigma} \Omega
= A^\sigma_{\nu} \ombar^\nu \wedge \Omega_\sigma = A^\sigma_{\nu} \delta_\sigma^\nu \, \Omega = (\Tr A_Z) \, \Omega \, . \tag*{\qed}
\end{equation*}
\end{pf*}

\begin{lemma}
The natural characterising form $\Omega$ evolves under the flow of $Z_v$ according to the formula
\begin{equation*}
\Lie_{Z_v} \Omega = v^i \pd{Z^\sigma_i}{y^\sigma} \Omega \, .
\end{equation*}
\end{lemma}
\begin{pf*}
From
\begin{equation*}
d \ombar^\sigma = - dZ^\sigma_i \wedge dx^i = - Z_k(Z^\sigma_i) dx^k \wedge dx^i - \pd{Z^\sigma_i}{y^\nu} \ombar^\nu \wedge dx^i =
\pd{Z^\sigma_i}{y^\nu} dx^i \wedge \ombar^\nu
\end{equation*}
we see that
\begin{align*}
d\Omega & = d (\ombar^1 \wedge \ombar^2 \wedge \cdots \wedge \ombar^m) \\
& = \sum_{\sigma=1}^m (-1)^{\sigma-1} \, \ombar^1 \wedge \cdots \wedge \ombar^{\sigma-1} \wedge d \ombar^\sigma \wedge \ombar^{\sigma+1} \wedge \cdots \wedge \ombar^m \\
& = \sum_{\sigma=1}^m (-1)^{\sigma-1} \pd{Z^\sigma_i}{y^\sigma} \ombar^1 \wedge \cdots \wedge \ombar^{\sigma-1} \wedge dx^i \wedge \ombar^\sigma \wedge \ombar^{\sigma+1} \wedge \cdots \wedge \ombar^m \\
& = \pd{Z^\sigma_i}{y^\sigma} dx^i \wedge \Omega \, .
\end{align*}
As $i_{Z_v} \ombar^\sigma = 0$ we see that $i_{Z_v} \Omega = 0$, so that
\begin{equation*}
\Lie_{Z_v} \Omega = i_{Z_v} d \Omega = v^i \pd{Z^\sigma_i}{y^\sigma} \Omega \, . \tag*{\qed}
\end{equation*}
\end{pf*}

\begin{prop}
The characterising $m$-form $\tilde{\Omega}$ is a $\Gamma$-compatible volume if, and only if, $\tilde{\Omega} = \mu \Omega$ where $\mu$ is a nowhere vanishing function such that
\begin{equation*}
Z_v(\mu) + \mu v^i \Bar H^\sigma_{\sigma i} = 0 \, .
\end{equation*}
\end{prop}
\begin{pf*}
Every characterising $m$-form for $D_Z$ is a nonzero multiple of $\Omega$, so that $\tilde{\Omega} = \mu \Omega$ for some nowhere vanishing function $\mu$; the compatibility condition then gives
\begin{align*}
\Lie_{Z_v} \tilde{\Omega} & = \Lie_{Z_v} (\mu \Omega)
= Z_v(\mu) \, \Omega + \mu \, \Lie_{Z_v} \Omega = v^i \biggl( Z_i(\mu) + \mu \pd{Z^\sigma_i}{y^\sigma} \biggr) \Omega \, , \\
i_{A_Z} \tilde \Omega & = i_{A_Z} (\mu \Omega) = \mu \, i_{A_Z} \Omega = \mu (\Tr A_Z) \,  \Omega \, .
\end{align*}
Using the formulas
\begin{equation*}
A^\sigma_\nu = v^i \biggl( \pd{Z^\sigma_i}{y^\nu} - \Hbar^\sigma_{\nu i} \biggr) \, , \qquad
\Tr A_Z = v^i \biggl( \pd{Z^\sigma_i}{y^\sigma} - \Hbar^\sigma_{\sigma i} \biggr)
\end{equation*}
we obtain
\begin{equation*}
0 = v^i \biggl( Z_i(\mu) + \mu \pd{Z^\sigma_i}{y^\sigma} \biggr) - \mu \Tr A_Z
= v^i \biggl( Z_i(\mu) + \mu \Hbar^\sigma_{\sigma i} \biggr) \, . \tag*{\qed}
\end{equation*}
\end{pf*}

\begin{theorem} \emph{(Second Main Theorem)}
\label{TSMT}

If $\tilde{\Omega}$ is a $\Gamma$-compatible volume then
\begin{equation*}
\Lie_{Z_v} \tilde{\Omega} = \Tr A_Z \, \tilde{\Omega} \, .
\end{equation*}
\end{theorem}
\begin{pf*}
A straightforward calculation gives
\begin{equation*}
\Lie_{Z_v} \tilde{\Omega} = Z_v (\mu) \Omega + v^i \pd{Z^\sigma_i}{y^\sigma} \tilde{\Omega}
= - v^i \Hbar^\sigma_{\sigma i} \tilde{\Omega} + v^i \pd{Z^\sigma_i}{y^\sigma} \tilde{\Omega}
= v^i \biggl(\pd{Z^\sigma_i}{y^\sigma} - \Hbar^\sigma_{\sigma i} \biggr) \tilde{\Omega}
= A^\sigma_\sigma \, \tilde{\Omega} \, . \tag*{\qed}
\end{equation*}
\end{pf*}

We see that this result is an exact analogue of the corresponding result for ODEs~\cite[Theorem~5.2]{JePr2000}, where the vector field $v$ on $X$ `slices' the congruence of solutions of $\Gamma$ into a congruence of curves.

\begin{theorem}
\label{Tcollapse}
Let $\gamma$ be an integral curve of $Z_v$ defined on some neighbourhood $I$ of zero, and let $C : I \to \R$ be the function $C(t) = (\Tr A_Z)(t)$. Fix a frame at $\gamma(0) \in X$ and let $\tilde{\Omega}(s)$ be the value of $\tilde{\Omega}$ on that frame Lie-transported to $\gamma(s)$ along $\gamma$. Then
\begin{equation}
\label{Edrag}
\Lie_{Z_v} \bigl( \tilde{\Omega}(t) \bigr) = C(t) \tilde{\Omega}(t) \, ,
\end{equation}
and if the integral
\begin{equation*}
\int_0^t C(s) ds
\end{equation*}
exists for $t \in [0, t_0)$ but diverges to $-\infty$ as $t \to t_0$ then $\tilde{\Omega}(t_0) = 0$.
\end{theorem}
\begin{pf}
Equation~\eqref{Edrag} is an immediate consequence of the previous theorem. As a formal solution to this equation is given by
\begin{equation*}
\tilde{\Omega}(t) = \tilde{\Omega}(0) \exp \biggl( \int_0^t C(s) ds \biggr)
\end{equation*}
the divergence of the integral corresponds to the vanishing of $\tilde{\Omega}$.
\end{pf}

As the vanishing of the $\Gamma$-compatible volume $\tilde{\Omega}$ indicates a collapse in the congruence of solution curves of the second order vector field $\Gamma_v$ specified by the vector field $Z_v$, it is a sufficient condition for the collapse in the congruence of solutions of the second order connection $\Gamma$ specified by the connection $Z$.


\subsection{The total shape operator}

The directional shape operator $A_Z$ is a vector-valued $1$-form n $Y$ which, from its coordinate expression
\begin{equation*}
A_Z = v^i \psibar^\sigma_i \otimes \vf{y^\sigma} \, ,
\end{equation*}
would appear to be the contraction of a vector-valued $2$-form with the vector field $v$. This is indeed the case, and we shall now see how to construct such a vector-valued $2$-form $\Ahat_Z$. As before,  we shall assume that the non-vanishing closed form $\varphi$ and the transversal vector field $v$ are fixed.

When defining the directional shape operator we used the map $Z^\varphi$ acting on vector-valued forms taking their values in $\dist_+$. Our strategy now will be to introduce a new operator similar to $Z^\varphi$, but acting on vector-valued forms taking their values in the larger bundle $V\pi_{1,0}$. As the dimension of $V\pi_{1,0}$ is greater than that of $V\pi$, there is no isomorphism between these bundles. We can, however, make use of the isomorphism
\begin{equation*}
\shat : \vf{y^\sigma_i} \mapsto dx^i \otimes \vf{y^\sigma} \, .
\end{equation*}
If $\Kform \in \Lambda^q(J^1\pi) \otimes V\pi_{1,0}$ and $\xi_1, \xi_2, \ldots, \xi_q, \xi_{q+1}$ are vector fields on $Y$ then, for any $y \in Y$, put
\begin{align*}
\lefteqn{(\Zupw \Kform)(y) \, (\xi_1, \xi_2, \ldots, \xi_q, \xi_{q+1})} \\
& \qquad := \sum_{r=1}^{q+1} (-1)^{q+r} i_{\xi_r} \biggl( \shat_{Z(y)} \Bigl( \Kform \bigl( Z(y) \bigr) \, \bigl( TZ(\xi_1), \ldots, \widehat{TZ(\xi_r)}, \ldots, TZ (\xi_{q+1}) \bigr) \Bigr) \biggr)
\end{align*}
where the wide circumflex denotes omission, so that $\Zupw \Kform \in  \Lambda^{q+1}(Y) \otimes V\pi$. If
\begin{equation*}
\Kform = \eta^\sigma_i \otimes \vf{y^\sigma_i} \, ,
\end{equation*}
then
\begin{equation*}
\Zupw \Kform = (-1)^{q+1}  dx^i \wedge Z^* \eta^\sigma_i \otimes \vf{y^\sigma} \, .
\end{equation*}

\begin{defn}
The vector valued $2$-form $\Ahat_Z$ on $Y$ defined by
\begin{equation*}
\Ahat_Z := \Zupw \vertproj
\end{equation*}
will be called a \emph{total shape operator} for the congruence defined by $Z$. \qed
\end{defn}
Note that, although $\Zupw$ does not depend on $\varphi$, the total shape operator $\Ahat_Z$ depends on both $\varphi$ and $v$, because $\vertproj$ does. In coordinates
\begin{equation*}
\Ahat_Z = Z^*(dx^i \wedge \psi^\sigma) \otimes \vf{y^\sigma} = dx^i \wedge \psibar^\sigma_i \otimes \vf{y^\sigma}
= \biggl( \pd{Z^\sigma_i}{y^\nu} - \Hbar_{\nu i}^\sigma \biggr) dx^i \wedge \ombar^\nu \otimes \vf{y^\sigma} \, ,
\end{equation*}
so that writing $\Ahat_Z = A^\sigma_{\nu i} dx^i \wedge \ombar^\nu \otimes \p / \p y^\sigma$ the components of $\Ahat_Z$ are given by
\begin{equation*}
A^\sigma_{\nu i} = \pd{Z^\sigma_i}{y^\nu} - \Hbar_{\nu i}^\sigma \, .
\end{equation*}

The relationship between the total shape operator and the directional shape operators defined earlier is straightforward.
\begin{prop}
If $Z$ is a congruence embedded in a second order connection $\Gamma$ then
\begin{equation*}
A_Z = i_{Z_v} \Ahat_Z \, .
\end{equation*}
\end{prop}
\begin{pf*}
As $i_{Z_v} \ombar^\sigma = 0$ (so that $i_{Z_v} \psibar^\sigma_i = 0$) and $i_{Z_v} dx^i = v^i$ we see immediately that
\begin{equation*}
i_{Z_v} \Ahat_Z = v^i \psibar^\sigma_i \otimes \vf{y^\sigma} = A_Z \, . \tag*{\qed}
\end{equation*}
\end{pf*}

Rewriting this formula in terms of the corresponding projectors, we obtain:
\begin{cor}
\begin{equation*}
Z^\varphi \vertproj_+ = i_{Z_v} \Zupw \vertproj \, .  \tag*{\qed}
\end{equation*}
\end{cor}

We now consider the analogues of Theorems~\ref{thmQ+} and~\ref{thmA} in the context of the total shape operator. The analogue of Theorem~\ref{thmQ+}, on the directional deformation of the vertical projector $\vertproj_+$, is Theorem~\ref{thmQ} on vertical curvature, and this may be interpreted as the deformation of the projector $\vertproj$ along $\Gamma$. We may then see that the analogue of Theorem~\ref{thmA} takes the following form, where $\nrb{\cdot,\cdot}$ denotes the algebraic \NR\ bracket of vector valued forms.

\begin{theorem}  \label{thmQA}  \emph{(Third Main Theorem)}

The total shape operator $\Ahat_Z$ satisfies the equation
\begin{equation*}
\fnb{\mathcal{Z}, \Ahat_Z} = - A^2 - \Zupw \Phi \,,
\end{equation*}
where $A^2$ denotes the vector-valued $2$-form
\begin{equation*}
A^2 = \tfrac{1}{4} \nrb{A_Z, A_Z} = \tfrac{1}{2} (A^\sigma_{\rho k} A_{\sigma i}^\nu - A^\sigma_{\rho i} A_{\sigma k}^\nu)
dx^i \wedge dx^k \wedge \ombar^\rho \otimes \vf{y^\nu} \, .
\end{equation*}
\end{theorem}

\begin{pf*}
We carry out the calculation in coordinates. First,
\begin{align*}
\fnb{\EHcon, \Ahat_Z} & = \biggfnb{dx^i \otimes Z_i, \, A^\nu_{\sigma k} dx^k \wedge \ombar^\sigma \otimes \vf{y^\sigma}} \\
& = dx^i \wedge dx^k \wedge \biggl( A^\nu_{\sigma k} \pd{Z^\sigma_i}{y^\rho} - A^\sigma_{\rho k} \pd{Z^\nu_i}{y^\sigma}
+ Z_i(A^\nu_{\rho k}) \biggr) \ombar^\rho \otimes \vf{y^\nu} \\
& = dx^i \wedge dx^k \wedge \biggl( \pd{Z^\sigma_i}{y^\rho} \pd{Z^\nu_k}{y^\sigma}
- \Hbar^\nu_{\sigma k} \pd{Z^\sigma_i}{y^\rho} - \Hbar^\sigma_{\rho i} \pd{Z^\nu_k}{y^\sigma}
- Z_i(\Hbar^\nu_{\rho k}) \biggr) \ombar^\rho \otimes \vf{y^\nu} \, ;
\end{align*}
on the other hand,
\begin{equation*}
\Zupw \Phi = - \Phibar^\nu_{i \rho k} dx^k \wedge dx^i \wedge \ombar^\rho \otimes \vf{y^\nu}
\end{equation*}
and
\begin{equation*}
\Phibar^\nu_{i \rho k} = \Gamma_i(H^\nu_{\rho k} ) \circ Z - H_\rho(F^\nu_{ik}) \circ Z + \Hbar^\sigma_{\rho i} \Hbar_{\sigma k}^\nu \, ,
\end{equation*}
so that, using $F_{ik} = F_{ki}$, we obtain
\begin{equation*}
\Zupw \Phi
= \bigl( Z_i(\Hbar^\nu_{\rho k}) + \Hbar^\sigma_{\rho i} \Hbar_{\sigma k}^\nu \bigr) dx^i \wedge dx^k \wedge \ombar^\rho \otimes \vf{y^\nu} \, .
\end{equation*}
Adding then gives
\begin{align*}
\fnb{\EHcon, \Ahat_Z} + \Zupw \Phi & = dx^i \wedge dx^k \wedge \biggl( \pd{Z^\sigma_i}{y^\rho} \pd{Z^\nu_k}{y^\sigma}
- \Hbar^\nu_{\sigma k} \pd{Z^\sigma_i}{y^\rho} - \Hbar^\sigma_{\rho i} \pd{Z^\nu_k}{y^\sigma}
+ \Hbar^\sigma_{\rho i} \Hbar_{\sigma k}^\nu \biggr) \ombar^\rho \otimes \vf{y^\nu} \\
& = A^\sigma_{\rho i} A_{\sigma k}^\nu \, dx^i \wedge dx^k \wedge \ombar^\rho \otimes \vf{y^\nu} \, . \tag*{\qed}
\end{align*}
\end{pf*}

Finally we give a result corresponding to Lemma~\ref{Lshapevol} for directional shape operators.
\begin{lemma}
Let
\begin{equation*}
\Omega = \ombar^1 \wedge \ombar^2 \wedge \cdots \wedge \ombar^m
\end{equation*}
be the natural characterising form for the distribution $D_Z$ in the coordinates $(x^i,y^\sigma)$. Then
\begin{equation*}
i_{\Ahat_Z} \Omega =  \Tr \Ahat_Z \wedge \Omega \, ,
\end{equation*}
where  $\Tr \Ahat_Z$, the `trace' of the total shape operator, is the $1$-form
\begin{equation*}
\Tr \Ahat_Z = A^\sigma_{\sigma i} \, dx^i \, .
\end{equation*}
\end{lemma}

\begin{pf*}
A straightforward calculation gives
\begin{equation*}
i_{\Ahat_Z} \Omega =  A^\sigma_{\nu i} dx^i \wedge \ombar^\nu \wedge i_{\p / \p y^\sigma} \Omega
= A^\sigma_{\nu i} dx^i \wedge \ombar^\nu \wedge \Omega_\sigma
= A^\sigma_{\nu i} \delta_\sigma^\nu dx^i \wedge \Omega =  \Tr \Ahat_Z \wedge \Omega \, . \tag*{\qed}
\end{equation*}
\end{pf*}

We also notice that there is a relationship between the traces of $A_Z$ and $\Ahat_Z$ given by
\begin{equation*}
\Tr A_Z = i_{Z_v} \Tr \Ahat_Z \, ,
\end{equation*}
so that
\begin{equation*}
i_{Z_v } i_{\Ahat_Z} \Omega =  (i_{Z_v} \Tr \Ahat_Z ) \, \Omega =  (\Tr A_Z) \, \Omega = i_{A_Z} \Omega \, .
\end{equation*}


\section{Applications and examples}

Although we have considered only PDE systems of connection type, the range of applications of our results is much broader. Given a single second order partial differential equation, or a system $\mathcal{E}$ of such second order PDEs, it is often easy to find, at least locally, a second order connection $\Gamma$ embedded in $\mathcal{E}$, so that every solution of $\Gamma$ is a solution of $\mathcal{E}$. For example, if $\mathcal{E}$ is a system of variational equations, then we know much more:\ as proved in \cite{RoSa14}, for every regular Lagrangian one has a family of second order connections, parametrised by traceless matrices, such that every integral section of each connection is an extremal, a solution of the Euler--Lagrange equations. Moreover, the integral sections of this family of second order connections include, at least locally, all extremals of the variational problem.

To give an elementary illustration how this approach of constructing an embedded second order connection can work, we revisit the simple lemniscate example from~\cite{SRP15} (see Figure~1).
\begin{figure}[ht]
\begin{center}
\includegraphics[width=20pc]{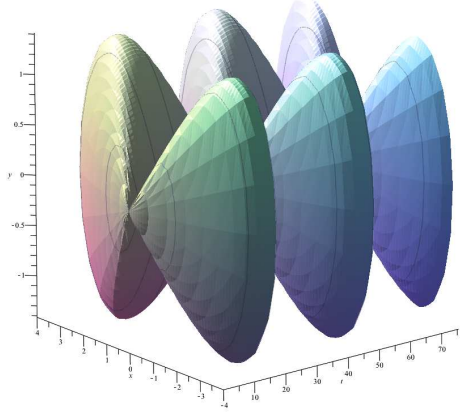}
\caption{Singularities of the lemniscate surface}\label{figure1}
\end{center}
\end{figure}
There, we considered a family of lemniscates with oscillating amplitude $a(t)$, given by
\begin{equation*}
r^2 = a^2(t) \cos 2\theta \, , \qquad \ddot{a} = -a \, ,
\end{equation*}
where $r$ was the dependent variable and $(t,\theta)$ were independent variables. We saw that, for each amplitude function $a$, this was a solution of a second order connection $\Gamma$ whose components were
\begin{equation*}
\Gamma_{\theta\theta} = - 2r - \frac{r_\theta^2}{r} \, , \qquad \Gamma_{t\theta} = \frac{r_t r_\theta}{r} \, , \qquad \Gamma_{tt} = -r \, ;
\end{equation*}
we regard this family of solutions as a congruence derived from a local connection $Z$, and use Theorem~\ref{Tcollapse} to analyse its collapse in the two coordinate directions.

We first compute the connection form $\EHcon = dt \otimes Z_t + d\theta \otimes Z_\theta$ using the two vector fields $Z_t$, $Z_\theta$ tangent to the family. We have $r_t \circ Z = r \cot t$ and $r_\theta \circ Z = -r \tan 2\theta$, giving
\begin{equation*}
Z_t = \vf{t} + Z^r_t \vf{r} = \vf{t} + r \cot t \vf{r} \, ,
\end{equation*}
and
\begin{equation*}
Z_\theta = \vf{\theta} + Z^r_\theta \vf{r} = \vf{\theta} - r \tan 2\theta \vf{r} \, .
\end{equation*}

We next construct directional shape operators, and to do this we need to make choices of $(\varphi, v)$. In~\cite{SRP15} we computed the (single) vector field $H_r$ spanning the distribution $\dist_-$ for two different choices:\ taking $(dt, \partial / \partial t)$ gave us
\begin{equation*}
H^{(t)}_r = \vf{r} + \frac{r_\theta}{r} \vf{r_\theta}
\end{equation*}
whereas taking $(d\theta, \partial / \partial \theta)$ gave us
\begin{equation*}
H^{(\theta)}_r = \vf{r} + \frac{r_t}{r} \vf{r_t} - \frac{r_\theta}{r} \vf{r_\theta} \, .
\end{equation*}
These choices give us the shape operators
\begin{equation*}
A^{(t)}_Z = \biggl( \pd{Z^r_t}{r} - Z^* \bigl( H^{(t)}_r(r_t) \bigr) \biggr) \ombar \otimes \vf{r} = (\cot t) \, \ombar \otimes \vf{r}
\end{equation*}
and
\begin{equation*}
A^{(\theta)}_Z = \biggl( \pd{Z^r_\theta}{r} - Z_\theta^* \bigl( H^{(\theta)}_r(r_\theta) \bigr) \biggr) \ombar \otimes \vf{r}
= (- 2 \tan 2\theta) \, \ombar \otimes \vf{r} \, .
\end{equation*}
Taking traces, we therefore obtain the functions
\begin{equation*}
C^{(t)}(s) = \bigl( \Tr A^{(t)}_Z \bigr)(s) = \cot s \, , \qquad C^{(\theta)}(s) = \bigl( \Tr A^{(\theta)}_Z \bigr)(s) = - 2 \tan 2s \, .
\end{equation*}
In the $t$ direction we have
\begin{equation*}
\int_{\pi/2}^t C^{(t)}(s) ds = \int_{\pi/2}^t \cot s \, ds = \int_{\pi/2}^t (\sin s)^{-1} \cos s \, ds = \log \sin t \, ;
\end{equation*}
on $[\pi/2, \pi)$ the integral is negative, but it diverges to $-\infty$ as $t \to \pi$. Thus for any angle $\theta$ the congruence of lemniscates collapses at $t = \pi$ (and at integer multiples of $\pi$).

In the $\theta$ direction we have
\begin{equation*}
\int_0^\theta C^{(\theta)}(s) ds = - 2 \int_0^\theta \tan 2s \, ds = \int_0^\theta (\cos 2s)^{-1} (- 2\sin s) \, ds
= \log \cos 2\theta \, ;
\end{equation*}
on $[0, \pi/4)$ the integral is zero or negative, but it diverges to $-\infty$ as $\theta \to \pi/4$. Thus at a fixed time $t$ the congruence collapses at $\theta = \pi/4$ (and at odd multiples thereof).


G.E.\,Prince \\
The Australian Mathematical Sciences Institute, c/o The University of Melbourne, Victoria 3010, Australia \\
\emph{and} Department of Mathematics and Statistics, La Trobe University, Victoria 3086, Australia \\
Email: \url{geoff@amsi.org.au}

O.\,Rossi \\
Department of Mathematics, Faculty of Science, University of Ostrava,\\ 30.\,dubna 22, 701 03 Ostrava, Czech Republic \\
\emph{and} Department of Mathematics and Statistics, La Trobe University, Victoria 3086, Australia \\
Email: \url{olga.rossi@osu.cz}

D.J.\,Saunders \\
Department of Mathematics, Faculty of Science, The University of Ostrava,\\ 30.\,dubna 22, 701 03 Ostrava, Czech Republic \\
Email: \url{david@symplectic.demon.co.uk}

\end{document}